\documentclass[11pt]{amsart}
\usepackage[margin=1in]{geometry}

\usepackage{amssymb}
\usepackage{amsthm}
\usepackage{amsmath}
\usepackage{mathrsfs}
\usepackage{amsbsy}
\usepackage[all]{xy}
\usepackage{bm}
\usepackage{hyperref}
\usepackage{tikz}
\usepackage{array}
\usepackage{float}
\usepackage{enumerate}
\usepackage{xcolor}
\usepackage{hhline}
\allowdisplaybreaks
\usepackage[noadjust]{cite}

\usepackage{caption}
\usepackage{tabu}
\usepackage{diagbox}
\usepackage{tikz}
\usepackage{comment}

\usepackage[noabbrev,capitalise]{cleveref}

\newenvironment{enumerate*}%
  {\begin{enumerate}[(I)]%
    \setlength{\itemsep}{10pt}%
    \setlength{\parskip}{0pt}}%
  {\end{enumerate}}

\newtheorem{theorem}{Theorem}[section]
\newtheorem{proposition}[theorem]{Proposition}

\newtheorem{conjecture}[theorem]{Conjecture}

\newtheorem{problem}[theorem]{Problem}
\newtheorem{lemma}[theorem]{Lemma}

\theoremstyle{definition}

\DeclareMathOperator{\ML}{ML}
\DeclareMathOperator{\vol}{vol}
\DeclareMathOperator{\covol}{covol}
\DeclareMathOperator{\PSL}{PSL}
\DeclareMathOperator{\PGL}{PGL}
\DeclareMathOperator{\acc}{acc}
\DeclareMathOperator{\den}{den}
\DeclareMathOperator{\mult}{mult}

\begin{document}

\title[]{The structure of Lonely Runner spectra}

\author[Vikram Giri]{Vikram Giri}
\address[]{Department of Mathematics, ETH Z\"urich, Z\"urich 8092, Switzerland}
\email{vikramaditya.giri@math.ethz.ch}

\author[Noah Kravitz]{Noah Kravitz}
\address[]{Department of Mathematics, Princeton University, Princeton, NJ 08540, USA}
\email{nkravitz@princeton.edu}

\maketitle

\begin{abstract}
For each closed subtorus $T$ of $(\mathbb{R}/\mathbb{Z})^n$, let $D(T)$ denote the (infimal) $L^\infty$-distance from $T$ to the point $(1/2,\ldots, 1/2)$.  The $n$-th Lonely Runner spectrum $\mathcal{S}(n)$ is defined to be the set of all values achieved by $D(T)$ as $T$ ranges over the $1$-dimensional subtori of $(\mathbb{R}/\mathbb{Z})^n$ that are not contained in the coordinate hyperplanes.  The Lonely Runner Conjecture predicts that $\mathcal{S}(n) \subseteq [0,1/2-1/(n+1)]$.  Rather than attack this conjecture directly, we study the qualitative structure of the sets $\mathcal{S}(n)$ via their accumulation points.  This project brings into the picture the analogues of $\mathcal{S}(n)$ where $1$-dimensional subtori are replaced by $k$-dimensional subtori or $k$-dimensional subgroups.
\end{abstract}

\section{Introduction}

\subsection{Lonely runners}

The Lonely Runner Problem of Wills \cite{wills} and Cusick \cite{cusick} is based on the following setup.  Suppose $n+1$ runners start at the same position on a unit-length circular track and then begin running around the track at pairwise distinct constant speeds.  We say that a runner is \emph{lonely} at a particular time if their distance from every other runner is at least $1/(n+1)$.  The Lonely Runner Conjecture predicts that each runner will be lonely at some time (allowing different times for different runners).

Consider the frame of reference of a single runner, and suppose that the $n$ other runners have speeds $v_1, \ldots, v_n$ relative to our fixed runner.  Define the \emph{maximum loneliness} for this set of speeds to be
$$\ML(v_1, \ldots, v_n):=\sup_{t \in \mathbb{R}} \min_{1 \leq i \leq n} \|tv_i\|_{\mathbb{R}/\mathbb{Z}},$$
where $\|x\|_{\mathbb{R}/\mathbb{Z}}$ denotes the distance from $x$ to the nearest integer.  Reformulated in this language, the Lonely Runner Conjecture asserts that
$$\ML(v_1, \ldots, v_n) \geq 1/(n+1)$$
for all nonzero real numbers $v_1, \ldots, v_n$.  Using Kronecker's Equidistribution Theorem, Bohman, Holzman, and Kleitman \cite{bohman} showed that it suffices to prove the Lonely Runner Conjecture for integer speeds; we will revisit this perspective below the fold.

The Lonely Runner Problem has attracted substantial attention since Wills first posed it in 1967.  Due to the work of many authors \cite{BW, CP, BGGST, bohman, Renault, BS}, the Lonely Runner Conjecture is known to hold for $n \leq 6$; unfortunately, none of the approaches to these small-$n$ cases seem to generalize to an arbitrary number of runners.  The trivial bound $\ML(v_1, \ldots, v_n) \geq 1/(2n)$ has also been improved slightly (see~\cite{tao} and the references therein), and the conjecture has been proven in the case where all of the runners have slow speeds (e.g., \cite{tao, bohman-peng, carl, ashwin-mehtaab}).  For other work on and around the Lonely Runner Problem, see the references in \cite{thesis}.

\subsection{Maximum loneliness spectra}

In this paper, we will be concerned with the perspective introduced in the second author's thesis \cite{thesis}.  The new motivating question was: What can we say about the set $\widetilde{\mathcal{S}}(n)$ of \emph{all} of the possible values of $\ML(v_1, \ldots, v_n)$?  For example, what is the second-smallest element after (conjecturally) $1/(n+1)$?  Can we describe all of $\widetilde{\mathcal{S}}(n)$ explicitly?  More qualitatively, what do the accumulation points of $\widetilde{\mathcal{S}}(n)$ look like?  Notice that each $\widetilde{\mathcal{S}}(n)$ is a subset of the interval $(0,1/2]$.  We will see later that $\widetilde{\mathcal{S}}(n)$ is a closed subset of the rationals.

The focus of the work \cite{thesis} was the following conjecture about the ``bottom'' part of the set $\widetilde{\mathcal{S}}(n)$.

\begin{conjecture}[Loneliness Spectrum Conjecture (\cite{thesis}, Conjecture 1.2)] \label{conj:spectrum}
For every natural number $n \geq 2$, we have
$$\widetilde{\mathcal{S}}(n) \cap (0,1/n)=\{s/(ns+1): s \in \mathbb{N}\}.$$
\end{conjecture}
One should think of the quantity $s/(ns+1)$ as the result of ``rounding down'' $1/n$ to the nearest multiple of $1/(ns+1)$.  An explicit construction from \cite{thesis} shows that all of the values $s/(ns+1)$ are contained in $\widetilde{\mathcal{S}}(n)$.  The conjecture has no content for $n=1$ and is easily verified for $n=2$.  One of the main results of \cite{thesis} is a proof of the conjecture for $n=3$.

Using computer experiments, Fan and Sun \cite{computer} discovered a family of examples that disproves Conjecture~\ref{conj:spectrum} for $n=4$.  In particular, they showed that
$$\ML(8, 4r+3, 4r+11, 4r+19)=\frac{2r+7}{4(2r+7)+2}$$
for every integer $r \geq 0$.  This provides a second infinite family of maximum loneliness values limiting to $1/n$ (at least for $n=4$) in a structured way.  Following the ``$1/n$ rounded down'' philosophy, Fan and Sun proposed the following natural weakening of Conjecture~\ref{conj:spectrum}.

\begin{conjecture}[Amended Loneliness Spectrum Conjecture (\cite{computer}, Conjecture 1.3]\label{conj:weak}
For every natural number $n \geq 2$, we have
$$\widetilde{\mathcal{S}}(n) \cap (0,1/n) \subseteq \{s/(ns+k): s \in \mathbb{N}, 1 \leq k \leq n\}.$$
\end{conjecture}

\subsection{Accumulation points}

On a qualitative level, both Conjecture~\ref{conj:spectrum} and its weakening Conjecture~\ref{conj:weak} predict that $1/n$ is the smallest accumulation point of $\widetilde{\mathcal{S}}(n)$.  It is no coincidence that this number $1/n$ coincides with the conjectural smallest value of $\widetilde{\mathcal{S}}(n-1)$.  Recall that, for a subset $\mathcal{S}$ of $\mathbb{R}$, the number $x \in \mathbb{R}$ is an \emph{upper accumulation point} if $\mathcal{S} \cap (x, x+\varepsilon) \neq \emptyset$ for all $\varepsilon>0$; we define \emph{lower accumulation points} analogously.  We write $\acc(\mathcal{S})$ for the set of (all) accumulation points of $\mathcal{S}$.  The following theorem forms the starting point for our investigations.

\begin{theorem}[\cite{thesis}, Theorem 6.5]\label{thm:accumulation-thesis}
Let $n \geq 2$ be a natural number.  Then the set of lower accumulation points of $\widetilde{\mathcal{S}}(n)$ contains $\widetilde{\mathcal{S}}(n-1)$.
\end{theorem}
Also in~\cite{thesis}, the second author raised the questions of whether upper accumulation points are impossible, and whether the set containment in Theorem~\ref{thm:accumulation-thesis} is actually an equality.  The goal of the present paper is to give an affirmative answer to the first part of this question and make partial progress towards the second part.  Along the way, our arguments will show that the sets $\widetilde{S}(n)$ and some closely related sets have a ``hierarchical'' structure; this information is potentially useful for inductive approaches to the Lonely Runner Conjecture, which so far have seemed inaccessible because of a failure to relate the instances of the conjecture for different numbers of runners.  Before stating our results precisely, we need to introduce a few geometric notions.

\subsection{View-obstruction, subtori, and subgroups}\label{sec:view}
It will be convenient for us to work with the ``view-obstruction'' formulation of the Lonely Runner Problem, as popularized by Cusick \cite{cusick}.  Associate each tuple $(v_1, \ldots, v_n)$ of nonzero integers with the $1$-dimensional subtorus $$T:=\pi(\langle (v_1, \ldots, v_n) \rangle_{\mathbb{R}})=\{t(v_1, \ldots, v_n): t \in \mathbb{R}\}/\mathbb{Z}^n$$ of the torus $(\mathbb{R}/\mathbb{Z})^n$, where $\pi: \mathbb{R}^n \to (\mathbb{R}/\mathbb{Z})^n$ is the standard quotient map.  Notice that $T$ is not contained in any of the coordinate hyperplanes of $(\mathbb{R}/\mathbb{Z})^n$; in general, we say that a subtorus is \emph{proper} if it is not contained in the union of the coordinate hyperplanes.  Let $D(T)$ denote the $L^\infty$-distance from $T$ to the point $(1/2, \ldots, 1/2)$ in the ``center'' of $(\mathbb{R}/\mathbb{Z})^n$, so that
$$D(T)=1/2-\ML(v_1, \ldots, v_n).$$
Now define the \emph{$n$-th Lonely Runner spectrum} to be the set $\mathcal{S}(n)$ of all values achieved by $D(T)$ as $T$ ranges over the $1$-dimensional proper subtori of $(\mathbb{R}/\mathbb{Z})^n$; we have
$$\mathcal{S}(n)=1/2-\widetilde{\mathcal{S}}(n).$$
In this new language, yet another reformulation of the Lonely Runner Conjecture is the assertion that $\mathcal{S}(n) \subseteq [0,1/2-1/(n+1)]$.  Notice that lower accumulation points for $\widetilde{\mathcal{S}}(n)$ correspond to upper accumulation points for $\mathcal{S}(n)$, and vice versa.

The more geometric line of inquiry arising from view-obstruction has been fruitfully related to the study of polytopes, and the ``zonotope'' version of the Lonely Runner Conjecture has received significant attention (see, e.g., \cite{HM, BHS, BNV, Beck-Schymura} and the references therein).

The function $D$ and the notion of properness also make sense for arbitrary closed subsets $X \subseteq (\mathbb{R}/\mathbb{Z})^n$.  In particular, define $D(X)$ to be the $L^\infty$-distance from $X$ to the point $(1/2, \ldots, 1/2)$, and say that $X$ is \emph{proper} if is is not contained in the union of the coordinate hyperplanes.\footnote{This condition is equivalent to $D(X)<1/2$, and it is slightly stronger than saying that $X$ is not contained in any single coordinate hyperplane.}  We will have occasion to study the values assumed by $D(X)$ not only when $X$ is a $1$-dimensional proper subtorus but also when $X$ is a higher-dimensional subtorus or a closed subgroup.

To this end, for $1 \leq k \leq n$, let $\mathcal{S}_k(n)$ be the set of all values achieved by $D(T)$ as $T$ ranges over the $k$-dimensional proper subtori of $(\mathbb{R}/\mathbb{Z})^n$.  Notice that $\mathcal{S}(n)=\mathcal{S}_1(n)$ in this notation.  We will also need to work with closed subgroups of $(\mathbb{R}/\mathbb{Z})^n$; in this paper, ``subgroup'' always means ``closed subgroup''.  Every subgroup $\Delta \subseteq (\mathbb{R}/\mathbb{Z})^n$ is the direct product of a $k$-dimensional subtorus of $(\mathbb{R}/\mathbb{Z})^n$ and a finite subgroup of $(\mathbb{R}/\mathbb{Z})^n$, and we say that the dimension of $\Delta$ is $k$.  Now, for $0 \leq k \leq n$, let $\mathcal{S}^*_k(n)$ denote the set of all values achieved by $D(\Delta)$ as $\Delta$ ranges over the $k$-dimensional proper subgroups of $(\mathbb{R}/\mathbb{Z})^n$.  We make the convention $\mathcal{S}^*_0(0):=\{0\}$, and we write $\mathcal{S}^*(n)$ for $\mathcal{S}^*_1(n)$.  Since every (proper) $k$-dimensional subtorus is also a (proper) $k$-dimensional subgroup, we have the inclusion $\mathcal{S}_k(n) \subseteq \mathcal{S}^*_k(n)$, which in general is strict (see below). 

Finally, we remark that one can think of $\mathcal{S}_k(n)$ as a multiset by keeping track of the values of $D(T)$ for all $k$-dimensional proper subtori $T \subseteq (\mathbb{R}/\mathbb{Z})^n$; we will denote this multiset by $\mathcal{S}_{k,\mult}(n)$.  The multiset $\mathcal{S}^*_{k,\mult}(n)$ is obtained from $\mathcal{S}^*_{k}(n)$ analogously.  The \emph{density points} of a multiset $\mathcal{S}$ are defined to be the accumulation points of $\mathcal{S}$ together with the infinite-multiplicity elements of $\mathcal{S}$.  We write $\den(\mathcal{S})$ for the set of density points of $\mathcal{S}$.

\subsection{Main result}
We are finally ready to state our main result, which is summarized by the following relationships among the sets $\mathcal{S}_k(n), \mathcal{S}^*_k(n)$.

\begin{theorem}\label{thm:relationships}
Let $n \geq 2$ be a natural number, and let $1 \leq k< n$ and $0 \leq k'< n$.  Then the sets $\mathcal{S}_k(n),\mathcal{S}_{k'}^*(n)$ have only upper accumulation points.  On the level of sets we have the inclusions
\begin{align*}
\acc(\mathcal{S}_k(n)) \subseteq \mathcal{S}_{k+1}(n) \subseteq \mathcal{S}_{k+1}^*(n)=\mathcal{S}_0^*(n-k-1),
\end{align*}
and on the level of multisets we have the equalities
$$\den(\mathcal{S}_{k,\mult}(n))=\mathcal{S}_{k+1}(n) \quad \text{and} \quad \den(\mathcal{S}_{k',\mult}^*(n))=\mathcal{S}_{k'+1}^*(n)=\mathcal{S}_{k'}^*(n-1).$$
Moreover, when $k=1$, we also have the inclusion $\mathcal{S}(n-1) \subseteq \acc(\mathcal{S}(n))$, which implies that
\begin{align}\label{eq:main-inclusions}
\mathcal{S}(n-1) \subseteq \acc(\mathcal{S}(n)) \subseteq \mathcal{S}_{2}(n) \subseteq \mathcal{S}_{2}^*(n)=\mathcal{S}^*(n-1).
\end{align}
\end{theorem}
The natural question, of course, is which of the inclusions  \eqref{eq:main-inclusions} are actually equalities.  Equality holds in \eqref{eq:main-inclusions} for $n=2$ trivially; it also holds for $n=3$ because (see Section~\ref{sec:examples} below) we can explicitly compute
$$\mathcal{S}(2)=\mathcal{S}^*(2)=\{0\} \cup \left\{1/(4s+2): s \in \mathbb{N}\right\}.$$
Equality cannot hold in general, however; we will show in Proposition~\ref{prop:0.14} that $7/50 \in \mathcal{S}^*(3) \setminus \mathcal{S}(3)$.  Nonetheless, we conjecture that equality holds in Theorem~\ref{thm:accumulation-thesis}.

\begin{conjecture}\label{conj:main-equality}
For any natural number $n \geq 2$, we have $\acc(\mathcal{S}(n))=\mathcal{S}(n-1)$.
\end{conjecture}

To prove this conjecture, it would suffice to show that $\mathcal{S}_2(n)=\mathcal{S}(n-1)$, i.e, that the first two inclusions in \eqref{eq:main-inclusions} are equalities.  More generally, we pose the following problem, which we suspect has an affirmative answer.

\begin{problem}\label{prob:higher-dim-subtori}
Is it the case that $\mathcal{S}_k(n)=\mathcal{S}_{k-\ell}(n-\ell)$ for all natural numbers $1\leq \ell<k \leq n$?
\end{problem}

Theorem~\ref{thm:relationships} has a couple of nice consequences.  One is that, due to the trivial containments $\mathcal{S}_k(n) \subseteq \mathcal{S}_{k-1}(n)$ and $\mathcal{S}^*_{k'}(n) \subseteq \mathcal{S}^*_{k'-1}(n)$, the sets $\mathcal{S}(n), \mathcal{S}^*(n)$ are closed for all $n$.  Another consequence is that each of $\mathcal{S}(n), \mathcal{S}^*(n)$ is a well-ordered set with order type $\omega^{n-1}+1$.

\subsection{Checking small speeds suffices}\label{sec:computability}


The maximum loneliness function is easy to compute, and it is natural to ask whether the Lonely Runner Conjecture can be verified by checking that there are no counterexamples with all speeds below some explicit threshold.  Tao \cite{tao} showed that this question can be answered in the affirmative, and that it suffices to check all speeds up to $n^{Cn^2}$ for some explicitly computable constant $C$.  Recently, Malikiosis, Santos, and Schymura \cite{MSS} used the zonotope formulation of the Lonely Runner Conjecture to improve this threshold to $(n(n+1)/2)^{n-1}\approx n^{2n}$.  The new perspective afforded by Theorem~\ref{thm:relationships} leads to a simple and transparent proof of the same result with the slightly weaker bound $n^{(5/2)n}$ (see below for our precise bound).  It is interesting to note that Minkowski's Theorem on successive minima plays a role both for us and for Malikiosis, Santos, and Schymura.

\subsection{Organization of the paper}

In Section~\ref{sec:similar}, we put the Lonely Runner spectra in their proper context alongside other ``bass note spectra'' such as the Markoff spectrum, and we sketch an interpretation in terms of abelian Bloch wave theory; this perspective was suggested to us by Peter Sarnak.  We prove our main result Theorem~\ref{thm:relationships} over the course of the following four sections: Section~\ref{sec:cubes} includes some geometric lemmas about cubes; Section~\ref{sec:kronecker} establishes a quantitative Kronecker-type result about high-volume subtori; Section~\ref{sec:producing} shows how to produce accumulation points in Lonely Runner spectra; and Section~\ref{sec:combining} assembles the pieces.  In Section~\ref{sec:small-speeds} we explain the computability result mentioned in Section~\ref{sec:computability}, and in Section~\ref{sec:examples} we provide several examples and characterizations in low-dimensional cases.  Finally, we raise some prospects for future work in Section~\ref{sec:remarks}.

\section{Relation to similar problems}\label{sec:similar}

\begin{figure}[b]
    \centering
    \begin{tikzpicture}
        \draw[<->] (-6,-2) -- (6,-2);
        \draw[<->] (-6,0) -- (6,0);
        \draw[<->] (-6,2) -- (6,2);
        \foreach \n in {1,...,100}
        {
        \filldraw[black] (7/\n-5,2) circle (1.2pt);
        }
        \foreach \n in {1,...,50}
        {
        \filldraw[black] (3/\n,0) circle (1.2pt);
        }
        \draw[line width=1mm , dotted, gray] (-2,0) -- (0,0);
        \draw[line width=1mm , black] (-5,0) -- (-2,0);
        \foreach \n in {1,...,100}
        {
        \filldraw[gray] (2/\n+2,-2) circle (1.0pt);
        }
        \foreach \n in {1,...,100}
        {
        \filldraw[gray] (2/\n-1.5,-2) circle (1.0pt);
        }
        \foreach \n in {1,...,100}
        {
        \filldraw[gray] (2/\n-2.66,-2) circle (1.0pt);
        }
        \foreach \n in {1,...,100}
        {
        \filldraw[gray] (2/\n-3.25,-2) circle (1.0pt);
        }
        \foreach \n in {1,...,100}
        {
        \filldraw[gray] (1/\n-3.6,-2) circle (1.0pt);
        }
        \foreach \n in {1,...,100}
        {
        \filldraw[black] (7/\n-5,-2) circle (1.5pt);
        }
        \filldraw[black] (-5,2) circle (2pt) node[anchor=north]{0};
        \filldraw[black] (-5,0) circle (2pt) node[anchor=north]{0};
        \filldraw[black] (-5,-2) circle (2pt) node[anchor=north]{0};
        \filldraw[black] (4,-2) circle (0pt) node[anchor=north]{1/4};
        \filldraw[black] (2,-2) circle (2pt) node[anchor=north]{1/6};
        \filldraw[black] (-1.5,-2) circle (2pt) node[anchor=north]{1/10};
        \filldraw[black] (-3.4,-2.13) circle (0pt) node[anchor=north]{...};
        \filldraw[black] (3,0) circle (0pt) node[anchor=north]{$1/\sqrt{5}$};
        \filldraw[black] (0,0) circle (0pt) node[anchor=north]{$1/3$};
        \filldraw[black] (-2,-0.05) circle (0pt) node[anchor=north]{$\mu_0$};
        \end{tikzpicture}
    \caption{These three sketches illustrate the behaviour of various bass-note spectra.  Top: The Lonely Runner spectrum $\mathcal{S}(2)$ is a ``rigid'' spectrum which is infinite and discrete with the unique accumulation point $0$.
    Middle: The Markoff spectrum has a continuous (``flexible'') bottom part, a discrete (``rigid'') top part, and a fractal transition part. Bottom: The Lonely Runner spectrum $\mathcal{S}(3)$ looks complicated, but it exhibits the ``hierarchical'' structure that its set of accumulation points is $\mathcal{S}(2)$.}
    \label{fig}
\end{figure}

\subsection{Abelian covers} \label{sec:abelian}
There is a natural generalization of the setup described in Section~\ref{sec:view}.  For any continuous function $f: (\mathbb{R}/\mathbb{Z})^n \to \mathbb{R}$, we can define $D_f(T)$ to be the minimum value of $f$ assumed on the subtorus $T \subseteq (\mathbb{R}/\mathbb{Z})^n$.  Our $D(T)$ from above corresponds to the function $f(x)=\Vert x-(1/2, \ldots, 1/2) \Vert_{\infty}$.  One can further define $\mathcal{S}_{f,k}(n)$ to be the set of all values assumed by $D_f(T)$ where $T$ is a $k$-dimensional proper subtorus of $(\mathbb{R}/\mathbb{Z})^n$, and it is interesting to study the qualitative structure of such sets.  (For some choices of $f$, it may make more sense to drop the ``properness'' condition.)

Abelian Bloch wave theory provides a fertile source of functions $f$.  Let $M$ be a geometric object such as a manifold or a graph, and let $\widetilde{M}$ be its universal abelian cover. Denote by $\pi_1(M)$ the fundamental group of $M$. The \emph{character torus} of $M$ is
$$T_M:=\{\text{characters $\chi: \pi_1(M) \to \mathbb{C}$ such that $\chi(\gamma)=1$ for all $\gamma \in \pi_1(\widetilde{M})$}\}$$
(which really is a torus in the usual sense of the term).  Finite abelian covers of $M$ correspond to finite subgroups of $T_M$, and infinite abelian covers correspond to closed subgroups of $T_M$ of dimension at least $1$.  The universal abelian cover of $M$ corresponds to the entire character torus.  The overarching principle of Bloch wave theory (or Floquet theory) is that spectral properties of an operator on an abelian cover $M'$ of $M$ can be understood in terms of ``twisted'' versions of the operator on the base space $M$; the subgroup of $T_M$ corresponding to $M'$ dictates which
character twists to consider.  In many cases of interest (such as when the operator is the Laplacian), these twisted operators are perturbations of the original operator, and Kato's perturbation theory (see \cite{kato}) tells us that quantities of interest vary continuously as we move along the character torus.  We can then define the function $f(\chi)$ to be the value of such a quantity for the twisted operator corresponding to the point $\chi$.  For example, taking $f(\chi)$ to be the spectral gap of the $\chi$-twisted Laplacian on $M$ results in the ``bass-note spectrum'' of $M$. See Section~$3$ of \cite{k-s} for an overview of how Bloch wave theory applies to the study of graph spectra; see also \cite{sspinor} and the work \cite{anshul} of Adve and the first author for applications to spinor spectra of closed Riemann surfaces; and see~\cite{girvin} for a standard treatment of Bloch wave theory from a physics perspective.

Our results and techniques in the context of the Lonely Runner Problem may also provide a path to understanding similar qualitative phenomena in spectra of differential operators on abelian covers.  Even if there turns out not to be a formal connection, it is helpful to situate the Lonely Runner Problem in this broader context.

\subsection{The Markoff spectrum}
In Section~\ref{sec:abelian} we described a question about (simple) geodesics in the flat torus $(\mathbb{R}/\mathbb{Z})^n$.  One can generalize this setup to geodesics in other locally symmetric spaces.

One particularly nice example comes from the locally symmetric space $X=\mathbb{H}^2/\Gamma$, where $\mathbb{H}^2:=\{x+iy: y>0\}$ is the complex upper half-plane with the hyperbolic metric and $\Gamma$ is a finite-index subgroup of the modular group $\PSL_2(\mathbb{Z})$ (which acts on $\mathbb{H}^2$ by M\"obius transformations).  Given a continuous function $f$ on $X$, we can define the function $D_{X,f}(\gamma)$ to be the minimum value of $f$ achieved on the closed (simple) geodesic $\gamma$, and one can ask about the set $\mathcal{S}_{X,f}$ of values of $D_{X,f}(\gamma)$.

The top part of the famous \emph{Markoff spectrum} (which is related to the Lagrange spectrum from Diophantine approximation) can be described in this language: It equals $\mathcal{S}_{X,f}$ for a certain choice of a congruence subgroup $\Gamma$ and a function $f$ for which $D_{X,f}(\gamma)$ measures how far $\gamma$ travels into the ``cusp'' of $X$ (i.e., the ``distance from $\gamma$ to the point at infinity'').  The structure of the Markoff spectrum is quite complicated (see Figure~\ref{fig}) and remains far from understood.  The top part of the spectrum is discrete and the bottom consists of a continuous interval (the so-called Hall--Freiman ray); there is a ``fractal'' region in the middle.  See~\cite{series} and Chapter~$7$ of~\cite{cusick-markoff} for overviews.

One can also consider the locally symmetric spaces $\PGL_n(\mathbb{R})/\PGL_n(\mathbb{Z})$ for $n\geq 3$. In many number-theoretic applications of this setup, it is natural to consider collections of geodesics (so-called ``packets'') or certain families of closed torus orbits rather than individual geodesics.  The \emph{measure rigidity conjectures} would imply that any infinite sequence of distinct such packets must eventually go arbitrarily deep into the cusp, and this would imply the rigidity of the corresponding spectra. See~\cite{e1, elmv1, elmv2} for more in this direction.

\section{Geometry of cubes}\label{sec:cubes}

We begin with some short ``combinatorial'' lemmas about how subtori and subgroups of $(\mathbb{R}/\mathbb{Z})^n$ can intersect cubes centered at $(1/2, \ldots, 1/2)$.  

\begin{lemma}\label{lem:rational}
Let $n \geq 2$ and $0 \leq k \leq n$ be natural numbers.  Then all elements of $\mathcal{S}^*_k(n)$ (and a fortiori of $\mathcal{S}_k(n)$ for $k \geq 1$) are rational numbers.  Moreover, for each $k$-dimensional proper subgroup $\Delta$ of $(\mathbb{R}/\mathbb{Z})^n$, there is a point $q \in \Delta$ with all rational coordinates such that $D(\Delta)=\Vert q-(1/2, \ldots, 1/2)\Vert_{\infty}$.
\end{lemma}

\begin{proof}
Let $\Delta$ be a $k$-dimensional proper subgroup of $(\mathbb{R}/\mathbb{Z})^n$.  Then there exist vectors $u_1, \ldots, u_k \in \mathbb{Z}^n$ and $v_1, \ldots, v_\ell \in \mathbb{Q}^n$ such that
$$\Delta=\bigcup_{1\leq j \leq \ell} \pi(\mathbb{R}u_1+\cdots+\mathbb{R}u_k + v_j)\,.$$
The point $p =\alpha_1 u_1+\cdots+\alpha_k u_k+v_j$ has image $\pi(p)=p-\lfloor p \rfloor$ (with $\lfloor p \rfloor$ taken coordinate-wise) in the fundamental domain $[0,1)^n$ of $(\mathbb{R}/\mathbb{Z})^n$, and hence the $L^\infty$-distance from $p$ to $(1/2, \ldots, 1/2)$ in $(\mathbb{R}/\mathbb{Z})^n$ is $$\Vert p-\lfloor p \rfloor-(1/2, \ldots, 1/2) \Vert_{L^\infty(\mathbb{R}^n)}.$$
Since the expression inside the $L^\infty$-norm is unchanged if the $\alpha_i$'s are shifted by integers (recall that the $u_i$'s have all integer coordinates), we have
$$D(\Delta)=\min_{1 \leq j \leq \ell}\min_{\alpha_1, \ldots, \alpha_k \in [0,1]} \left\Vert \sum_i \alpha_i u_i + v_j -\left\lfloor \sum_i \alpha_i u_i + v_j \right\rfloor-(1/2, \ldots, 1/2) \right\Vert_{L^{\infty} (\mathbb{R})^n}.$$
As $j$ and $\alpha_1, \ldots, \alpha_k$ range, the quantity $\lfloor \sum_i \alpha_i u_i + v_j \rfloor$ ranges over finitely many vectors in $\mathbb{Z}^n$, say, $m_1, \ldots, m_R \in \mathbb{Z}^n$.
Then we can write
$$D(\Delta)=\min_{1 \leq r \leq R}\min_{1\leq j \leq \ell} \min_{\alpha_1, \ldots, \alpha_k \in [0,1]} \left\Vert \sum_i \alpha_i u_i+v_j -m_r-(1/2, \ldots, 1/2) \right\Vert_{L^{\infty} (\mathbb{R})^n}.$$
The expression inside the $L^\infty$-norm is a rational affine-linear function of the $\alpha_i$'s for each choice of $r,j$, so we conclude that the set of $\alpha_i$'s attaining the minimum is a finite union of rational polytopes (i.e., polytopes defined by linear equalities and inequalities with rational coefficients).  Hence some minimizer has all rational coordinates, and it follows that $D(T)$ is also rational.
\end{proof}

It is clear that $\mathcal{S}_{k-\ell}^*(n-\ell) \subseteq \mathcal{S}^*_{k}(n)$ for all $1 \leq \ell \leq k \leq n$: If $\Delta$ is a $(k-\ell)$-dimensional proper subgroup of $(\mathbb{R}/\mathbb{Z})^{n-\ell}$ with $D(\Delta)=d$, then $\Delta':=\Delta \times (\mathbb{R}/\mathbb{Z})^\ell$ is a $k$-dimensional proper subgroup of $(\mathbb{R}/\mathbb{Z})^n$ with $D(\Delta')=d$.  Somewhat surprisingly, this inclusion is an equality.  It is this ``hierarchical'' structure that makes the function $D$ particularly nice (among all of the possible functions $D_f$ described in Section~\ref{sec:abelian}).  One should think of the equality $\mathcal{S}^*_{k}(n)=\mathcal{S}_{k-\ell}^*(n-\ell)$ as saying that the set of attainable values of $D(\Delta)$ depends on only the codimension of the subgroups $\Delta$ under consideration.

\begin{lemma}\label{lem:cubes}
For natural numbers $1 \leq \ell\leq  k\leq n$, we have $\mathcal{S}^*_{k}(n)=\mathcal{S}_{k-\ell}^*(n-\ell)$.
\end{lemma}

\begin{proof}
The result is trivial for $k=n$, so suppose that $k<n$.  By the previous remark, it suffices to show that $\mathcal{S}^*_{k}(n)\subseteq \mathcal{S}_{k-\ell}^*(n-\ell)$.  It is clear that $0 \in \mathcal{S}_{k-\ell}^*(n-\ell)$, as witnessed by, for instance, the subgroup $$\{(\underbrace{0, \ldots, 0}_{n-k}), (\underbrace{1/2, \ldots, 1/2}_{n-k})\} \times (\mathbb{R}/\mathbb{Z})^{k-\ell} \subseteq (\mathbb{R}/\mathbb{Z})^{n-\ell}.$$  Now let $d \in \mathcal{S}_{k}^*(n)$ with $d>0$, and let $\Delta$ be a $k$-dimensional proper subgroup of $(\mathbb{R}/\mathbb{Z})^n$ such that $D(\Delta)=d$.  Let $B$ denote the box in $(\mathbb{R}/\mathbb{Z})^n$ consisting of the points whose $L^\infty$-distance from $(1/2, \ldots, 1/2)$ is $d$.  The properness of $\Delta$ ensures that $d<1/2$, so $B$ really is a ``box''.  We claim that $\Delta$ intersects a face of $B$ of dimension at most $n-k-1$.  Indeed, let $F$ be a lowest-dimensional face of $B$ that intersects $\Delta$. If $\dim(F)>n-k$, then the intersection of $F$ and $\Delta$ contains a line segment by naive dimension-counting, and we can follow this line segment to find an intersection of $\Delta$ with a lower-dimensional face of $B$.  If $\dim(F)=n-k$, then again the intersection of $F$ and $\Delta$ contains a line segment, since otherwise $\Delta$ would intersect $F$ transversely and pass into the interior of $B$, which is impossible.  So we conclude that $\dim(F) \leq n-k-1$.

Let $p$ be an intersection point of $\Delta$ with a face of $B$ of dimension at most $n-k-1$.  Then $p$ has at least $k+1$ coordinates with values in $\{1/2+d,1/2-d\}$; without loss of generality, we may assume that the first $k+1 \geq \ell+1$ coordinates of $p$ equal $1/2+d$.  Recall from Lemma~\ref{lem:cubes} that $d$ is rational and that the set of points $x \in \Delta$ with $\Vert x-(1/2, \ldots, 1/2) \Vert_\infty=d$ is a finite union of rational polytopes.  If we impose the additional constraint that the first $\ell+1$ coordinates all equal $1/2+d$, then we again obtain a finite union of rational polytopes, and we know that it is nonempty since it contains $p$.  Hence it also contains some point $q$ with all rational coordinates.  Let $\Delta'$ be the subgroup of $\Delta$ obtained by intersecting $\Delta$ with the subspace $$\{(x_1,\ldots,x_n): x_1=\cdots=x_{\ell+1}\}.$$  Notice that $\Delta'$ is a subgroup of dimension at least $k-\ell$, and it is proper because it contains the point $q$.  Now, take any $(k-\ell)$-dimensional subgroup of $\Delta'$ containing $q$, and let $\Delta'' \subseteq (\mathbb{R}/\mathbb{Z})^{n-\ell}$ denote its the projection onto the last $n-\ell$ coordinates of $(\mathbb{R}/\mathbb{Z})^n$.  Then $\Delta'' \subseteq (\mathbb{R}/\mathbb{Z})^{n-\ell}$ is a $(k-\ell)$-dimensional proper subgroup with $D(\Delta'')=D(\Delta')=D(\Delta)=d$, so $d \in \mathcal{S}_{k-\ell}^*(n-\ell)$, as desired.
\end{proof}

Let us now describe what happens when one works with subtori instead of subgroups.  We still have the trivial inclusion $\mathcal{S}_{k-\ell}(n-\ell) \subseteq \mathcal{S}_k(n)$ for $1\leq \ell< k \leq n$, although one cannot hope for any relationship when $\ell=k$ since $\mathcal{S}_0(n-\ell)=\emptyset$, as there are no $0$-dimensional proper subtori.

Matters around the reverse inclusion are more complicated.  For $\ell<k$, one can run the first paragraph of the argument from the proof of Lemma~\ref{lem:cubes} on a $k$-dimensional proper subtorus $U$.  The trouble comes in the second paragraph, where the intersection of $U$ with the subspace $\{(x_1,\ldots,x_n): x_1=\cdots=x_{\ell+1}\}$ may be a \emph{disconnected} subgroup rather than a connected subtorus.\footnote{This subtlety was the cause of an error in an earlier version of this paper.}  Since $\ell<k$, one could also consider the intersections of $U$ with other subspaces of the form
$$\{(x_1, \ldots, x_n): x_{i_1}=\cdots=x_{i_{\ell+1}}\}$$
for $1 \leq i_1<\cdots<i_{\ell+1} \leq k+1$.  Unfortunately, even this additional flexibility is insufficient to guarantee that the intersection is connected; a concrete counterexample with $\ell=1, k=2, n=3$ is given by the subtorus $\pi(\langle (0,7,-5), (5,-93,0) \rangle_\mathbb{R}) \subseteq (\mathbb{R}/\mathbb{Z})^3$.

It is still possible to say a little bit. A special case of the following lemma appeared, in somewhat different language, as Lemma 8 of the paper of Bohman, Holzman, and Kleitman \cite{bohman}.  One should think of the lemma as as a characterization of $\max \mathcal{S}_k(n)$ by the ``codimension'' $n-k$.

\begin{lemma}\label{lem:BHK}
For natural numbers $1 \leq \ell< k \leq n$, we have $\max \mathcal{S}_k(n) =\max \mathcal{S}_{k-\ell}(n-\ell)$.
\end{lemma}

\begin{proof}
By induction, it suffices to establish the $\ell=1$ case of the lemma.  The trivial inclusion $\mathcal{S}_{k-1}(n-1) \subseteq \mathcal{S}_k(n)$ gives $\max \mathcal{S}_k(n) \geq \max \mathcal{S}_{k-1}(n-1)$.  It remains to show the reverse inequality.  

Let $d \in \mathcal{S}_k(n)$, and let $U \subseteq (\mathbb{R}/\mathbb{Z})^n$ be a $k$-dimensional proper subtorus with $D(U)=d$.  Write $U=\pi(\langle w_1,\ldots, w_{k-2},u,v \rangle_{\mathbb{R}})$ for nonzero vectors $w_1,\ldots, w_{k-2},u=(u_1, \ldots, u_n), v=(v_1, \ldots, v_n) \in \mathbb{Z}^n$.  The properness of $U$ guarantees that for each $1 \leq i \leq n$, at least one of the $k$ generators of $U$ has its $i$-th coordinate nonzero.  Thus, by replacing $u$ with $u$ plus suitable integer multiples of the other generators, we may assume that $u$ has all nonzero coordinates.  Moreover, by changing coordinates $x_i \mapsto -x_i$ for the $i$'s with $u_i<0$, we may assume that $u$ has all coordinates strictly positive.  Finally, by permuting the coordinates $x_i$, we may assume that the quantities
$$v_1/u_1, v_2/u_2, \ldots, v_n/u_n$$
are non-decreasing.  Since $u,v$ are non-parallel, there is some index $1 \leq i \leq n-1$ such that $v_i/u_i<v_{i+1}/u_{i+1}$.

Now, let $T$ denote the identity component of the intersection $U \cap \{x_i =-x_{i+1}\}$.  Note that $T$ is a subtorus of dimension either $k-1$ or $k$.  We claim that $T$ is proper.  To this end, the key observation is that $T$ contains the $1$-dimensional proper subtorus
$$\pi(\langle (v_i+v_{i+1})u-(u_i+u_{i+1})v \rangle_{\mathbb{R}}).$$
Indeed, the $i$-th and $(i+1)$-th coordinates of $(v_i+v_{i+1})u-(u_i+u_{i+1})v$ sum to zero by direct calculation, so this is indeed a subtorus of $T$.  For properness, it suffices to show that $(v_i+v_{i+1})u-(u_i+u_{i+1})v$ has all coordinates nonzero.  To see this, notice that the $j$-th coordinate is equal to
$$(v_i+v_{i+1})u_j-(u_i+u_{i+1})v_j,$$
which vanishes if and only if
$$\frac{v_j}{u_j}=\frac{v_i+v_{i+1}}{u_i+u_{i+1}}.$$
Since $u_i,u_{i+1}>0$ and $v_i/u_i<v_{i+1}/u_{i+1}$, the quantity $(v_i+v_{i+1})/(u_i+u_{i+1})$ lies strictly between $v_i/u_i$ and $v_{i+1}/u_{i+1}$, and in particular it is not equal to $v_j/u_j$ for any $1 \leq j \leq n$.

Take $T'$ to be a $(k-1)$-dimensional proper subtorus of $T$, and let $T''$ denote the projection of $T'$ onto all but the $i$-th coordinate.  Then $T''$ is a $(k-1)$-dimensional proper subtorus of $(\mathbb{R}/\mathbb{Z})^{n-1}$, and $D(U) \leq D(T)\leq D(T')=D(T'') \leq \max \mathcal{S}_{k-1}(n-1)$, as desired.
\end{proof}

We remark that we could have intersected $U$ with the subspace where $x_1=x_n$ instead of the subspace where $x_i=-x_{i+1}$; there is some (limited) flexibility in this choice.

It is instructive to compare what the proofs of Lemmas~\ref{lem:cubes} and~\ref{lem:BHK} tell us about the intersections $U':=U \cap \{x_i=x_j\}$ where  $U \subseteq (\mathbb{R}/\mathbb{Z})^n$ is a fixed $2$-dimensional proper subtorus and $1 \leq i<j \leq n$ range.  The proof of Lemma~\ref{lem:cubes} provides indices $i<j$ such that $U'$ satisfies $D(U')=D(U)$ but may be disconnected.  The proof of Lemma~\ref{lem:BHK} provides indices $i<j$ such that $U$ is connected but may satisfy $D(U')>D(U)$.  The example described between the two lemmas shows that for some choices of $U$ it is impossible to find indices $i<j$ such that $U'$ both is connected and satisfies $D(U')=D(U)$.

\section{Quantitative Kronecker Theorem}\label{sec:kronecker}

We will require some technical results on the distribution of subtori of $(\mathbb{R}/\mathbb{Z})^n$.  In this setting, every $k$-dimensional subtorus ($1 \leq k \leq n$) is of the form $T=\pi(W)$ for some $k$-dimensional subspace $W$ of $\mathbb{R}^n$ that is ``rational'' in the sense that $\Lambda:=\mathbb{Z}^n \cap W$ is a lattice of rank $k$.  We write $\vol_k(T)$ for the $k$-dimensional volume of a $k$-dimensional subtorus $T$.  Note that the volume of the torus $T=\pi(W)$ is simply the covolume of its associated lattice $\Lambda$. 

We will require a simple lemma showing that there are only finitely many subtori with volume smaller than any constant.  Schmidt obtained this result with a sharp bound in~\cite{schmidt}, but we will also include a short proof here (with a worse bound) to keep the paper self-contained.  We first recall a standard fact from reduction theory.
\begin{theorem}[\cite{siegel}, Sections X.5--6]\label{thm.red}
Let $\Lambda \subset \mathbb{R}^k$ be a lattice of full rank.  Then there is a basis $b_1, \ldots, b_k$ of $\Lambda$ such that
$$\Vert b_1\Vert_2 \cdots \Vert b_k \Vert_2 \leq 2^k (3/2)^{k(k-1)/2} \covol_k(\Lambda)/\omega_k,$$
where $\omega_k:=\pi^{k/2}/\Gamma(k/2+1)$ is the volume of the $L^2$-unit ball in $\mathbb{R}^k$.

\end{theorem}
This theorem is a variation of Minkowski's Second Theorem from the geometry of numbers.  We cannot apply Minkowski's Second Theorem directly because the elements of $\Lambda$ achieving the successive minima need not form a $\mathbb{Z}$-basis (as discussed in~\cite{siegel}, Section X.5).

\begin{lemma}\label{lem:finitesubtori}
Let $0 \leq k \leq n$ be nonnegative integers, and let $V>0$ be a positive real number.  Then there are only finitely many $k$-dimensional subtori of $(\mathbb{R}/\mathbb{Z})^n$ with volume at most $V$.
\end{lemma}

\begin{proof}

Every $k$-dimensional subtorus $T$ of $(\mathbb{R}/\mathbb{Z})^n$ with volume at most $V$ is of the form $\pi(W)$ for some $k$-dimensional subspace $W$ of $\mathbb{R}^n$ that is ``rational'' in the sense that $\Lambda:=\mathbb{Z}^n\cap W$ is a rank-$k$ lattice in $W$; identifying $W$ with $\mathbb{R}^k$, we see that $T$ is isomorphic to $\mathbb{R}^k/\Lambda$.  Theorem~\ref{thm.red} provides a basis $b_1, \ldots, b_k$ of $\Lambda$ such that
$$\Vert b_1\Vert_2 \cdots \Vert b_k \Vert_2 \leq 2^k (3/2)^{k(k-1)/2} \vol_k(T)/\omega_k.$$
Since every nonzero element of $\mathbb{Z}^n$ has length at least $1$, we see that each
$$\Vert b_i \Vert_2 \leq 2^k (3/2)^{k(k-1)/2}V/\omega_k:=\ell(k,V).$$
Since $\mathbb{Z}^n$ has only finitely many elements of length at most $\ell(k,V)$, we conclude that there are only finitely many choices for the $b_i$'s and hence (since $T$ is determined by $W$, which is determined by $\Lambda$) only finitely many choices for $T$.
\end{proof}

We also need the analogous result for subgroups.

\begin{lemma}\label{lem:finitesubgroups}
Let $0 \leq k \leq n$ be nonnegative integers, and let {$V>0$} be a positive real number.  Then there are only finitely many $k$-dimensional subgroups of $(\mathbb{R}/\mathbb{Z})^n$ of volume at most $V$.
\end{lemma}

\begin{proof}
Every subgroup $\Delta$ of $(\mathbb{R}/\mathbb{Z})^n$ with volume at most $V$ can be written as a direct sum $\Delta= T \oplus H$, where $T$ is the identity component of $\Delta$ and $H$ is a finite subgroup of $(\mathbb{R}/\mathbb{Z})^n$; note that $T$ is a $k$-dimensional subtorus and $H \subseteq (\mathbb{Q}/\mathbb{Z})^n$.  From $\vol_k(T) \geq 1$ and $$\vol_k(T) \cdot |H|=\vol_k(\Delta)\leq V,$$ we see that $\vol_k(T), |H| \leq V$.  Lemma~\ref{lem:finitesubgroups} tells us that there are only finitely many choices for $T$.  To see that there are only finitely many choice for $H$, note that $H$ is contained in the set of elements of $(\mathbb{R}/\mathbb{Z})^n$ of order at most $V$ and this set is finite.  Since $\Delta$ is determined by $T,H$, we conclude that there are only finitely many choices for $\Delta$.
\end{proof}

We can now prove our Kronecker-type result.  The classical version of Kronecker's Equidistribution Theorem says that every irrational orbit in $(\mathbb{R}/\mathbb{Z})^n$ equidistributes in some subtorus of dimension at least $2$.  Our version can be understood as a quantitative, finitary analogue for closed subtori and subgroups.

\begin{lemma}\label{lem:kronecker''}
Let $0 \leq k \leq n$ be nonnegative integers, and let $\varepsilon>0$ be a positive real number.  Then there exist a constant $C^* = C^*(n,k,\varepsilon) >0$ and a finite list $L^*= L^*(n,k,\varepsilon)$ of subgroups of $(\mathbb{R}/\mathbb{Z})^n$ such that for each $k$-dimensional subgroup $\Delta$ of $(\mathbb{R}/\mathbb{Z})^n$, one of the following holds:
\begin{enumerate}
    \item $\Delta$ has volume at most $C$.
    \item $\Delta$ is contained in one of the subgroups $\Gamma \in L^*$ and is $\varepsilon$-dense (with respect to the $L^2((\mathbb{R}/\mathbb{Z})^n)$-norm) in $\Gamma$. If $\Delta$ is a subtorus, then $\Gamma$ can also be taken to be subtorus.
\end{enumerate}
Moreover, each subgroup in $L^*$ has dimension at least $k+1$ if $k<n$.
\end{lemma}
Notice that $(\mathbb{R}/\mathbb{Z})^n$ will always be one of the elements of our list $L^*$; this corresponds to the situation where $\Delta$ is dense in $(\mathbb{R}/\mathbb{Z})^n$.

\begin{proof}
We fix $n$ and proceed by downward induction on $k$.  For the base case $k=n$, the only subgroup $\Delta$ is $\Delta=(\mathbb{R}/\mathbb{Z})^n$; take $L^*(n,n,\varepsilon)=\{(\mathbb{R}/\mathbb{Z})^n\}$ and $C^*(n,n,\varepsilon)=1$ (say), and there is nothing to prove.

We proceed to the induction step.  For visualizing the argument, the reader may find it helpful to keep the $k=n-1$ case in mind.  Let $\Delta$ be a $k$-dimensional subgroup of volume $V$, and let $T$ be the identity component of $\Delta$.  Note that $T$ is itself a $k$-dimensional subtorus, and write $T=\pi(W)$ with $W$ a $k$-dimensional rational subspace of $\mathbb{R}^n$.  Moreover, we can write $\Delta= T \oplus H$ for some finite group $H \subseteq (\mathbb{R}/\mathbb{Z})^n$, and there is a unique finite set $R \subseteq [0,1)^n$ such that $\pi(R)=H$ bijectively; note that in fact $R \subseteq \mathbb{Q}^n$.  Let $X$ denote the tubular neighborhood of $\Delta$ of radius $\varepsilon/2$.  If
$$V\omega_{n-k} (\varepsilon/2)^{n-k}> 1,$$
then the $(n-k)$-dimensional ``orthogonal slices'' of $X$ around the points of $\Delta$ cannot all be disjoint.  Lifting this picture to $\mathbb{R}^n$ gives us some $x \in (\mathbb{Z}^n+R) \setminus W$ such that the $\varepsilon/2$-tubular neighborhoods around $W$ and $W+x$ intersect.  In particular, the orthogonal complement of $W$ intersects $W+x$ at some (necessarily nonzero) point $p$ of length at most $\varepsilon$.  The point $p$ has all rational coordinates since it is determined by a rational system of linear equations (involving $W$ and $x$).  Let $W'$ be the $\mathbb{R}$-subspace of $\mathbb{R}^n$ spanned by $W$ and $p$, and let $T':=\pi(W')$.  Notice that $T'$ is a $(k+1)$-dimensional subtorus (because of the rationality of $p$) and that the $\varepsilon/2$-neighborhood of $\Delta$ contains $T'$.  It follows that $\Delta$ is $\varepsilon/2$-dense in the subgroup $\Delta':=T'+H$.

We now have a dichotomy depending on whether the volume of $\Delta'$ is small or large.  Set
$$L^*(n,k,\varepsilon):=L^*(n,k+1, \varepsilon/2) \cup \{\text{$(k+1)$-dim'l subgroups of volume at most $C^*(n,k+1,\varepsilon/2)$}\}.$$
If the volume of $\Delta'$ is at most $C^*(n,k+1,\varepsilon/2)$, then $\Delta' \in L^*(n,k,\varepsilon)$ and we are done since $\Delta$ is $\varepsilon$-dense in $\Delta'$.  If instead the volume of $\Delta'$ is larger than $C^*(n,k+1,\varepsilon/2)$, then there is some $\Delta'' \in L^*(n,k+1, \varepsilon/2) \subseteq L^*(n,k,\varepsilon)$ such that $\Delta'$ is $\varepsilon/2$-dense in $\Delta''$, and again we are done since $\Delta$ is $\varepsilon$-dense in $\Delta''$ by the Triangle Inequality.  This proves the lemma with the choice $$C^*(n,k,\varepsilon):=\frac{1}{\omega_{n-k} (\varepsilon/2)^{n-k}}.$$
Lemma~\ref{lem:finitesubgroups} implies that $L^*(n,k,\varepsilon)$ is a finite set. The second part of (2) is clear from tracing through the proof with only subtori under consideration.
\end{proof}

\section{Producing accumulation points}\label{sec:producing}
Our next task is showing how elements of certain spectra lead to accumulation points in other spectra.  We begin by producing density points in multiset spectra.

\begin{proposition}\label{prop:acc-pts-from-higher-dim}
For natural numbers $1 \leq k<n$, we have $\mathcal{S}_{k+1}(n) \subseteq \den(\mathcal{S}_{k,\mult}(n))$.  For nonnegative integers $0 \leq k< n$, we have $\mathcal{S}^*_{k+1}(n) \subseteq \den(\mathcal{S}^*_{k,\mult}(n))$. 
\end{proposition}

\begin{proof}
For the first statement, we must show that for every $(k+1)$-dimensional proper subtorus $U \subseteq (\mathbb{R}/\mathbb{Z})^n$, the value $D(U)$ is a density point of the multiset $\mathcal{S}_{k,\mult}(n)$.  Fix such a $U$.  Notice that every $k$-dimensional subtorus $T \subseteq U$ has $D(T) \geq D(U)$ by the definition of $D$.  If moreover $T$ is $\varepsilon$-dense in $U$ with respect to the $L^\infty$-norm, then the Triangle Inequality also gives $D(T) \leq D(U)+\varepsilon$.  Hence, it suffices to exhibit an infinite sequence of $k$-dimensional proper subtori $T_1, T_2, \ldots$ that are contained in $U$ and become $o(1)$-dense in $U$.  Notice that the properness condition comes for free if $T_1, T_2, \ldots$ become $o(1)$-dense in $U$, since non-proper subtori all have $D$-value $1/2$.

Now, in fact, any infinite sequence of $k$-dimensional subtori of $U$ becomes $o(1)$-dense in $U$ (by the codimension-$1$ cases of Lemmas~\ref{lem:finitesubtori} and~\ref{lem:kronecker''}).  We can also produce such a sequence explicitly by writing $U=\pi(\langle u_1, \ldots, u_{k+1} \rangle_{\mathbb{R}})$ for some nonzero vectors $u_1, \ldots, u_{k+1} \in \mathbb{Z}^n$ and setting $$T_j:=\pi(\langle u_1, u_2, \ldots, u_{k-1}, u_k+ju_{k+1} \rangle_{\mathbb{R}}) \subseteq U$$
for each $j \in \mathbb{N}$.  The $o(1)$-denseness of the $T_j$'s in $U$ follows from the $o(1)$-denseness of the subtori $\pi(\langle (1,j) \rangle_{\mathbb{R}})$ in $(\mathbb{R}/\mathbb{Z})^2$, which is geometrically obvious.


The second statement goes in the same way.  Let $\Delta \subseteq (\mathbb{R}/\mathbb{Z})^n$ be a $(k+1)$-dimensional proper subgroup.  As in the previous paragraph, it suffices to find an infinite sequence $\Gamma_1, \Gamma_2, \ldots$ of $k$-dimensional subgroups of $\Delta$ that become $o(1)$-dense in $\Gamma$.  We can write $\Delta$ as a direct sum $\Delta=U \oplus G$, where $U$ is a $(k+1)$-dimensional subtorus (not necessarily proper) and $G$ is a finite subgroup (also not necessarily proper).  
It suffices to produce a sequence of $k$-dimensional subgroups $H_1, H_2, \ldots$ of $U$ that become $o(1)$-dense in $U$, since then $\Gamma_1:=H_1 \oplus G, \Gamma_2:=H_2 \oplus G, \ldots$ will be $o(1)$-dense in $\Delta$.  If $k=0$, then $U$ is a $1$-dimensional subtorus and we can take $H_j$ to be the unique subgroup of $U$ of order $j$.  If $k \geq 1$, then we can take $H_j$ to be the subtorus $T_j$ constructed in the previous paragraph.
\end{proof}

The situation becomes more subtle when we work with set spectra instead of multiset spectra.  The proof strategy for Proposition~\ref{prop:acc-pts-from-higher-dim} completely breaks down.  In the context of the first statement (about subtori), it is possible that $D(T)=D(U)$ for all but finitely many $k$-dimensional proper subtori $T$ of $U$; in this case, the $D$-values of the $k$-dimensional proper subtori of $U$ witness $D(U)$ as a density point of $\mathcal{S}_{k,\mult}(n)$ but do not witness $D(U)$ as a genuine accumulation point of $\mathcal{S}_k(n)$.  Jain and the second author~\cite{vanshika} have given a simple geometric criterion for when a $2$-dimensional proper subtorus $U \subseteq (\mathbb{R}/\mathbb{Z})^n$ has infinitely many $1$-dimensional proper subtori $T$ with $D(T)>D(U)$: This occurs if and only if the locus where $U$ attains its $D$-value is contained in a finite union of parallel line segments.\footnote{The same argument shows that a $(k+1)$-dimensional proper subtorus $U$ has infinitely many $k$-dimensional proper subtori $T$ with $D(T)>D(U)$ if and only if the locus where $U$ attains its $D$-value is contained in a finite union of parallel $k$-dimensional disks.}  Using this criterion, they showed that, for example, the $2$-dimensional proper subtorus $$U=\pi(\langle (1, 2, 3, 2, 0, 0, 0),(0, 0, 0, 2, 1, 2, 3) \rangle_{\mathbb{R}}) \subseteq (\mathbb{R}/\mathbb{Z})^7$$ has only finitely many $1$-dimensional proper subtori $T$ with $D(T)>D(U)$.  It is even easier to construct such examples in the context of the second statement of Proposition~\ref{prop:acc-pts-from-higher-dim} (about subgroups).  For instance, the $1$-dimensional proper subgroup
$$\Delta=\mathbb{R}/\mathbb{Z} \times \{0,1/3,2/3\} \subseteq (\mathbb{R}/\mathbb{Z})^2$$
has only finitely many discrete subgroups $\Gamma$ with $D(\Gamma)>D(\Delta)$.

An important special case occurs when $U$ is of the form $U=U' \times \mathbb{R}/\mathbb{Z}$ for $U'$ a $1$-dimensional proper subtorus.  Since $U'$ achieves its $D$-value at only finitely many points, we see that $U$ achieves its $D$-value on only finitely many ``vertical'' line segements.  The criterion of Jain and second author guarantees that $U$ contains infinitely many $1$-dimensional proper subtori $T$ with $D(T)>D(U)=D(U')$.  In particular, there exist such $T$'s with volume tending to infinity, in which case $D(T)$ approaches $D(U)$ from above due to the codimension-$1$ cases of Lemmas~\ref{lem:finitesubtori} and~\ref{lem:kronecker''}.  We conclude that, as observed more directly in \cite{thesis}, we have the inclusion
$$\mathcal{S}(n-1) \subseteq \acc (\mathcal{S}(n))$$
for every natural number $n \geq 2$.  It is an interesting open problem to determine whether or not 
$$\mathcal{S}_k(n-1) \subseteq \acc (\mathcal{S}_k(n))$$
for all $1 \leq k<n$; we suspect that this is likely the case, and that in fact equality holds.

\section{Putting everything together}\label{sec:combining}

We have nearly all of the pieces that comprise Theorem~\ref{thm:relationships}.  Before we can complete the proof, we need to use our Kronecker-type result to relate $\acc (\mathcal{S}_k(n))$ to $\mathcal{S}_{k+1}(n)$.

\begin{lemma}\label{lem:hard-half}
Let $1 \leq k<n$ be natural numbers.  Then the accumulation points of $\mathcal{S}_k(n)$ are all accumulation points only from above, and $\den(\mathcal{S}_{k,\mult}(n)) \subseteq \mathcal{S}_{k+1}(n)$.  
\end{lemma}

\begin{proof}
We first claim that if $1 \leq \ell \leq n-1$ and $T_1, T_2, \ldots$ is a sequence of distinct $\ell$-dimensional proper subtori of $(\mathbb{R}/\mathbb{Z})^n$, then there exist a subtorus $T$ of dimension at least $\ell+1$ and a subsequence $T_{i_1}, T_{i_2}, \ldots$ such that each $T_{i_j}$ is contained in $T$ and the $T_{i_j}$'s become $o(1)$-dense in $T$ as $j \to \infty$.  We proceed by downward induction on $\ell$.  Notice that the volume of the $T_i$'s tends to infinity as $i$ grows.  The base case $\ell=n-1$ now follows immediately from Lemma~\ref{lem:kronecker''}, which tells us that the $T_i$'s are $o(1)$-dense in $(\mathbb{R}/\mathbb{Z})^n$.  For the induction step, suppose that the $T_i$'s are $\ell$-dimensional with $\ell<n-1$.  Lemma~\ref{lem:kronecker''} tells us that there are subtori $T'_i$, each of dimension at least $\ell+1$, such that each $T_i$ is contained in $T'_i$ and the $T_i$'s are $o(1)$-dense in their respective $T'_i$'s as $i \to \infty$.  By passing to a subsequence, we may assume that all of the $T'_i$'s have the same dimension, say, $\ell'\geq \ell+1$.  If there are only finitely many distinct $T'_i$'s, then one of them appears infinitely often and we can conclude by letting $T$ be such a $T'_i$.  If instead there are infinitely many distinct $T'_i$'s, then (by the induction hypothesis) there exist a subtorus $T$ of dimension at least $\ell'+1$ and a subsequence $T'_{i_1}, T'_{i_2}, \ldots$ such that each $T'_{i_j}$ is contained in $T$ and the $T'_{i_j}$'s are $o(1)$-dense in $T$.  Then we also see that each $T_{i_j}$ is contained in $T$ and the $T_{i_j}$'s are $o(1)$-dense in $T$ (by the Triangle Inequality), as desired.

We now prove the lemma by (upward) induction on $k$.  Let $T_1, T_2, \ldots$ be a sequence of $k$-dimensional proper subtori with $\lim_{i \to \infty}D(T_i)=d$.  The claim from the previous paragraph tells us that, after passing to a subsequence, we may assume that all of the $T_i$'s are contained in a single subtorus $T$ of dimension at least $k+1$ and that the $T_i$'s are $o(1)$-dense in $T$.  It follows that every $D(T_i)\geq D(T)$ and that $\lim_{i \to \infty} D(T_i)=D(T)$, so $D(T)=d$.  In particular, $D(T_i)$ cannot limit to $d$ from below.  This establishes the first statement of the lemma, namely, that $\mathcal{S}_k(n)$ has only upper accumulation points. 

We know that $T$ has dimension at least $k+1$.  We can find a $(k+1)$-dimensional subtorus $T'$ such that $d=D(T)=D(T') \in \mathcal{S}_{k+1}(n)$: Just take any $(k+1)$-dimensional subtorus $T'$ of $T$ that passes through a rational point of $T$ with $L^\infty$-distance $d$ from the point $(1/2, \ldots, 1/2)$ (using Lemma~\ref{lem:rational}). This establishes the second statement of the lemma, namely, that $\acc(\mathcal{S}_k(n)) \subseteq \mathcal{S}_{k+1}(n)$.
\end{proof}

The following subgroup-version of this lemma can be proven by an identical argument.

\begin{lemma}\label{lem:hard-half-subgroups}
Let $0 \leq k<n$ be nonnegative integers.  Then the accumulation points of $\mathcal{S}_k^*(n)$ are all accumulation points only from above, and $\den(\mathcal{S}_{k,\mult}^*(n)) \subseteq \mathcal{S}_{k+1}^*(n)$.
\end{lemma}

We remark that Lemmas~\ref{lem:hard-half} and~\ref{lem:hard-half-subgroups} do not use anything special about the fact that $D(T)$ measures the $L^\infty$ distance from $T$ to $(1/2, \ldots, 1/2)$: We would get the same ``qualitative'' result for any $D_f(T)$.

Combining all of our results so far yields Theorem~\ref{thm:relationships}, as follows.
\begin{itemize}
    \item Lemmas~\ref{lem:hard-half} and~\ref{lem:hard-half-subgroups} tell us that $\mathcal{S}_k(n), \mathcal{S}^*_{k'}(n)$ have only upper accumulation points.
    \item For the long string of set inclusions in Theorem~\ref{thm:relationships}:  We obtain $\acc (\mathcal{S}_k(n)) \subseteq \mathcal{S}_{k+1}(n)$ from Lemma~\ref{lem:hard-half}.  The inclusion $\mathcal{S}_{k+1}(n) \subseteq \mathcal{S}^*_{k+1}(n)$ is trivial.  And we obtain $\mathcal{S}^*_{k+1}(n)=\mathcal{S}_0^*(n-k-1)$ from Lemma~\ref{lem:cubes}.
    \item For the multiset equalities: We obtain $\mathcal{S}_{k+1}(n) \subseteq \den(\mathcal{S}_{k,\mult}(n))$ and $\mathcal{S}^*_{k'+1}(n) \subseteq \den(\mathcal{S}^*_{k',\mult}(n))$ from Proposition~\ref{prop:acc-pts-from-higher-dim}.  And we obtain $\den(\mathcal{S}_{k,\mult}(n)) \subseteq \mathcal{S}_{k+1}(n)$ and $\den(\mathcal{S}^*_{k',\mult}(n)) \subseteq \mathcal{S}^*_{k'+1}(n)$ from Lemmas~\ref{lem:hard-half} and~\ref{lem:hard-half-subgroups}.  The final equality $\mathcal{S}_{k'+1}^*(n)=\mathcal{S}_{k'}^*(n-1)$ follows from Lemma~\ref{lem:cubes}.
    \item For the final assertion, we obtain $\mathcal{S}(n-1) \subseteq \acc (\mathcal{S}(n))$ from Theorem~\ref{thm:accumulation-thesis}, and Lemma~\ref{lem:cubes} gives us $\mathcal{S}^*_2(n)=\mathcal{S}^*(n-1)$.
\end{itemize}

\section{Checking small speeds suffices}\label{sec:small-speeds}
In this short section, we explain how to deduce the computability result described in Section~\ref{sec:computability}.  The argument is a straightforward modification of the approach taken in the proof of Lemma~\ref{lem:kronecker''}.  Suppose we already know that the Lonely Runner Conjecture holds for $n-1$ runners and we want to determine whether or not it holds for $n$ runners, simply by checking the maximum loneliness for all $n$-tuples of small natural-number speeds.  It would suffice to show that if there exists a counterexample to the Lonely Runner Conjecture with $n$ runners, then there exists such a counterexample with all slow runners.  We will show a stronger statement, namely, that every $n$-tuple with large sum of squares of speeds must satisfy the Lonely Runner Conjecture.  Recall that the tuple of speeds $(v_1, \ldots, v_n)$ corresponds to the $1$-dimensional proper subtorus $T:=\pi(\langle(v_1, \ldots, v_n)\rangle_{\mathbb{R}}) \subseteq (\mathbb{R}/\mathbb{Z})^n$, which has length $\vol_1(T)=\sqrt{v_1^2+\cdots+v_n^2}$ if $\gcd(v_1, \ldots, v_n)=1$ (which of course we may assume).  

Let $\varepsilon:=1/n-1/(n+1)=1/(n(n+1))$, and let $T$ be a $1$-dimensional proper subtorus of $(\mathbb{R}/\mathbb{Z})^n$ with volume $V$.  The argument from the proof of Lemma~\ref{lem:kronecker''} shows that if
$$V \omega_{n-1}\varepsilon^{n-1}>1,$$
then there is some $2$-dimensional proper subtorus $T'$ such that $T$ is $\varepsilon$-dense in $T'$.  In particular, $D(T) \leq D(T')+\varepsilon$, and Lemma~\ref{lem:BHK} (together with our ``induction'' hypothesis) ensures that the latter quantity is at most $1/2-1/(n+1)$.  Hence $T$ satisfies the Lonely Runner Conjecture once
$$V>\frac{1}{\omega_{n-1} \varepsilon^{n-1}}=\frac{\Gamma((n+1)/2)(n(n+1))^{n-1}}{\pi^{(n-1)/2}},$$
and Stirling's Approximation shows that this threshold is smaller than $n^{(5/2)n}$.

We remark that we can run the same argument if the Lonely Runner Conjecture for $n-1$ runners fails but we still have $\max \mathcal{S}(n-1)=\max \mathcal{S}_2(n) \leq  1/2-1/(n+1)-\delta$ for some $\delta>0$; the threshold for $V$ will just be corresponding larger.

The ideas in this section also bear on the problem of characterizing \emph{tight instances} of the Lonely Runner Problem, i.e., $1$-dimensional proper subtori $T \subseteq (\mathbb{R}/\mathbb{Z})^n$ with the (conjecturally) largest possible value $D(T)=1/2-1/(n+1)$.  The above proof shows that if the Lonely Runner Conjecture holds for $n-1$ runners, then there are only finitely many tight instances for $n$ runners, and one could easily extract an explicit upper bound on their volume (and, by extension, their number).

\section{Calculations for small (co)dimension}\label{sec:examples}
When the number of runners is small, we can provide some more precise information about the relation between Lonely Runner spectra and subgroup Lonely Runner spectra; this leads to more general observations about (subgroup) Lonely Runner spectra of small codimension.

We start with $\mathcal{S}_1(1)=\{0\}$.  It is also easy to compute (see, e.g., \cite[Theorem 4.1]{thesis} or the remark on page 5 of \cite{BGGST}) that
$$\mathcal{S}_1(2)=\{0\} \cup \{1/(4s+2): s \in \mathbb{N}\}.$$
In particular,
$$\acc (\mathcal{S}_1(2))=\{0\}=\mathcal{S}_1(1),$$
which provides an affirmative answer to the $k=n$ (``codimension-$0$'') case of Problem~\ref{prob:higher-dim-subtori}.

Next, we calculate that
$$\mathcal{S}_0^*(1)=\{0\} \cup \{1/(4s+2): s \in \mathbb{N}\}=\mathcal{S}_1(2).$$
Indeed, each proper discrete subgroup of $\mathbb{R}/\mathbb{Z}$ is equal to $\pi(\langle 1/q \rangle_\mathbb{Z})$ for some natural number $q \geq 2$.  If $q$ is even, then $1/2 \in \pi(\langle 1/q \rangle_\mathbb{Z})$ and $D(\pi(\langle 1/q \rangle_\mathbb{Z}))=0$.  If instead $q=2s+1$ is odd, then the closest points in $\pi(\langle 1/q \rangle_\mathbb{Z})$ to $1/2$ are $\pm s/q$, so $$D(\pi(\langle 1/q \rangle_\mathbb{Z}))=1/2-s/q=1/(4s+2).$$
Notice that the degenerate Lonely Runner spectrum $\mathcal{S}_0(1)=\emptyset$ is emphatically not equal to $\mathcal{S}_0^*(1)$.

For each $n \geq 2$ we have the chain of inclusions
$$\mathcal{S}_1(2) \subseteq \mathcal{S}_{n-1}(n) \subseteq \mathcal{S}_{n-1}^*(n)=\mathcal{S}_0^*(1)=\mathcal{S}_1(2),$$
so equality must hold.  This proves the $k=n-1$ (``codimension-$1$'') case of Problem~\ref{prob:higher-dim-subtori} and implies that $\acc (\mathcal{S}_1(3))=\mathcal{S}_1(2)$, as expected.

The subgroup Lonely Runner spectrum $\mathcal{S}_0^*(2)$ is more complicated since there are many discrete subgroups of $(\mathbb{R}/\mathbb{Z})^2$.  Of course, we have the inclusion
$$\mathcal{S}_1(3) \subseteq \mathcal{S}_1^*(3)=\mathcal{S}_0^*(2).$$
This inclusion, however, is strict.
\begin{proposition}\label{prop:0.14}
We have $7/50 \in \mathcal{S}_0^*(2) \setminus \mathcal{S}_1(3)$.
\end{proposition}

\begin{proof}
A direct calculation shows that $D(\pi(\langle(12/25,9/25) \rangle_\mathbb{Z}))=7/50 \in \mathcal{S}_0^*(2)$; see Figure~\ref{fig:0.14-example}.  It remains to show that $7/50 \notin \mathcal{S}_1(3)$.  Let $T \subseteq (\mathbb{R}/\mathbb{Z})^3$ be a $1$-dimensional proper subtorus.  We will show that $D(T) \neq 7/50$.  First, suppose that $\vol(T)>199$; note that $199\cdot ((1/25-\varepsilon)^{2}\pi)>1$ for some $\varepsilon>0$.  Then the tubular neighborhood argument from the proof of Lemma~\ref{lem:kronecker''}, with a tube of radius $1/25-\varepsilon$, shows that $T$ is $(1/25-\varepsilon)$-dense in some proper subtorus $U \subseteq (\mathbb{R}/\mathbb{Z})^3$ of dimension at least $2$.  Thus
$$D(U) \leq D(T) \leq D(U)+1/25-\varepsilon.$$
We know from the characterization of $\mathcal{S}_2(3)$ that either $D(U)=1/6$ or $D(U) \leq 1/10$.  In the former case, the inequality $D(T) \geq D(U)=1/6$ implies that $D(T) \neq 7/50$.  In the latter case, the inequality
$$D(T) \leq D(U)+1/25-\varepsilon \leq 1/10+1/25-\varepsilon=7/50-\varepsilon$$ implies that $D(T) \neq 7/50$.  Finally, an exhaustive computer search (taking about an hour on a standard laptop) shows that $D(T) \neq 7/50$ when $\vol(T) \leq 199$.
\end{proof}

\begin{figure}[h]
    \centering
\begin{tikzpicture}[scale=3]
    \draw (0,0) -- (1,0) -- (1,1) -- (0,1) -- cycle;
    
    \fill[blue!30] 
        (0.36, 0.36) -- 
        (0.64, 0.36) -- 
        (0.64, 0.64) -- 
        (0.36, 0.64) -- cycle;
    
    \draw 
        (0.36, 0.36) -- 
        (0.64, 0.36) -- 
        (0.64, 0.64) -- 
        (0.36, 0.64) -- cycle;

    \foreach \x/\y in {0.48/0.36, 0.96/0.72, 0.44/0.08, 0.92/0.44, 
                        0.4/0.8, 0.88/0.16, 0.36/0.52, 0.84/0.88, 
                        0.32/0.24, 0.8/0.6, 0.28/0.96, 0.76/0.32, 
                        0.24/0.68, 0.72/0.04, 0.2/0.4, 0.68/0.76, 
                        0.16/0.12, 0.64/0.48, 0.12/0.84, 0.6/0.2, 
                        0.08/0.56, 0.56/0.92, 0.04/0.28, 0.52/0.64, 0/0}
    \fill[black] (\x, \y) circle (0.65pt);
\end{tikzpicture}
    \caption{The black dots show elements of the discrete subgroup $\pi(\langle (12/25,9/25) \rangle_\mathbb{Z}) \subseteq (\mathbb{R}/\mathbb{Z})^2$, which witnesses the value $7/50 \in \mathcal{S}_0^*(2)$, and the blue square shows the $L^\infty$-ball of radius $7/50$ centered at the point $(1/2,1/2)$.  Notice that there are black dots on all four edges of the blue square; this is a necessary feature of any such example.}\label{fig:0.14-example}
\end{figure}

We believe that $7/50 \notin \mathcal{S}_2(4)$, but we do not have the computing power to carry out an exhaustive computer search analogous to what we did in Proposition~\ref{prop:0.14}.  In any case, the example in Proposition~\ref{prop:0.14} shows that proving $\mathcal{S}_2(4)=\mathcal{S}_1(3)$ will necessarily be more difficult than proving $\mathcal{S}_2(3)=\mathcal{S}_1(2)$.  In particular, one would have to show that not all $1$-dimensional proper subgroups of $(\mathbb{R}/\mathbb{Z})^3$ can arise when one applies the slicing argument of Lemma~\ref{lem:cubes} to $2$-dimensional proper subtori in $(\mathbb{R}/\mathbb{Z})^4$.

\section{Remarks and open problems}\label{sec:remarks}
We conclude with a few comments and questions.
\begin{enumerate}
    \item Our proof of Theorem~\ref{thm:relationships} shows that each element of $\mathcal{S}_1(n) \setminus \mathcal{S}_2(n)$ appears with finite multiplicity (but not uniformly bounded multiplicity, even for $n=2$) in $\mathcal{S}_{1,\mult}(n)$.  In particular, as we mentioned in Section~\ref{sec:small-speeds}, if the Lonely Runner Conjecture is true for $n$ and $n-1$ runners, then there are only finitely many tight instances of the Lonely Runner Conjecture for $n$ runners, and all such instances are low-volume.  It would feasible (and very useful) to carry out numerical calculations in this direction.
    \item For a proper subtorus $U \subseteq (\mathbb{R}/\mathbb{Z})^n$ of dimension at least $2$, one can define the \emph{Lonely Runner spectrum relative to $U$} to be the set $\mathcal{S}(U)$ of all values of $D(T)$ as $T$ ranges over $1$-dimensional proper subtori contained in $U$.  Our proof of Theorem~\ref{thm:relationships} shows that all accumulation points of $\mathcal{S}(n)$ are ``due'' to relatively Lonely Runner spectra, in the sense that for each $d \in \acc(\mathcal{S}(n))$ there are some small $\varepsilon>0$ and a finite list $U_1, \ldots, U_t$ with $D(U_1)=\cdots=D(U_t)=d$  such that $$\mathcal{S}(n) \cap (d,d+\varepsilon)=\left(\bigcup_{j=1}^t \mathcal{S}(U_j) \right) \cap (d,d+\varepsilon).$$
    Jain and the second author \cite{vanshika} have studied such relatively Lonely Runner spectra for $2$-dimensional subtori $U$ and shown that the sets $\mathcal{S}(U)$ have surprisingly rigid arithmetical properties.
    \item Our results show that $$\acc (\mathcal{S}(3))=\mathcal{S}(2)=\{1/(4s+2): s \in \mathbb{N}\} \cup \{0\}.$$ The main result of~\cite{thesis} explicitly determines $\mathcal{S}(3) \cap [1/6, 1/2]$, i.e., the set $\mathcal{S}(3)$ ``up to the first accumulation point $1/6$'' (see Conjecture~\ref{conj:spectrum}). The follow-up work of Jain and the second author \cite{vanshika} has yielded a description of $\mathcal{S}(3)$ up to the second accumulation point $1/10$ as well as chunks of some other Lonely Runner spectra.  It would be interesting to obtain more such characterizations.
    \item Is the set $\cup_{n \geq 1}\mathcal{S}(n)$ dense in the interval $[0,1/2]$?  If not, does $\overline{\cup_{n \geq 1}\mathcal{S}(n)}$ have any ``nice'' self-similarity properties? 
\end{enumerate}

\section*{Acknowledgments}
We are grateful to Peter Sarnak for suggesting the use of a Kronecker-type lemma and for explaining the perspective that led to Section~\ref{sec:similar} of this paper.  We thank Alex Fan and Alec Sun, and Romanos Diogenes Malikiosis for sharing preprints with us.  We thank Terry Tao for drawing our attention to the work leading to \cite{MSS}, and we thank Giorgos Kotsovolis for explaining to us the measure rigidity conjectures.  We thank Milan Haiman for pointing out a gap in an earlier version of the proof of Lemma~\ref{lem:hard-half}, and we thank Vanshika Jain for helpful conversations.  The first author was supported by the National Science Foundation under grant DMS-FRG-1854344. The second author was supported in part by an NSF Graduate Research Fellowship (grant DGE--2039656).

\end{document}